\documentclass[letterpaper, 10pt, conference]{IEEEtran}

\IEEEoverridecommandlockouts                              
\makeatletter
	\let\NAT@parse\undefined
\makeatother

\PassOptionsToPackage{noend}{algpseudocode}
\PassOptionsToPackage{nameinlink,capitalize}{cleveref}
\usepackage[cal=boondoxo,scr=euler]{mathalfa}
\usepackage{booktabs}
\usepackage[space]{cite}
\usepackage{csvsimple}
\let\labelindent\relax
\usepackage[%
	thmreset=never,
	eqreset=never,
	alglinenumber,
	noload={caption,subfig,subcaption,todonotes},
]{myPreamble}
\usepackage{pgfplotstable}
\ifams
	\usepackage[foot]{amsaddr}
\fi

\newcommand{\affil}[1]{\texorpdfstring{\textsuperscript{#1}}{}}

	\allowdisplaybreaks

	\DeclareMathAlphabet{\mathcalligra}{T1}{calligra}{m}{n}
	\makeatletter
		\def\operator@font{\rm}
	\makeatother
	
	\newcommand{\msmall}[1]{\text{\small\(#1\)}}
	\theoremstyle{plain}
	\CreateNewTheorem{problem}{Problem}
	
	\newtheorem*{problem*}{Main problem}

	\makeatletter
		\renewcommand{\clevethm@proofsectiontitle}{Omitted proofs of }
		\let\@old@cref\cref
		\let\@old@Cref\Cref
		\renewcommand{\cref}[1]{%
			\IfBeginWith{#1}{proof:}{%
				the \hyperref[#1]{Appendix}%
			}{%
				\@old@cref{#1}%
			}%
		}
		\renewcommand{\Cref}[1]{%
			\IfBeginWith{#1}{proof:}{%
				The \hyperref[#1]{Appendix}%
			}{%
				\@old@Cref{#1}%
			}%
		}
	\makeatother
	\crefname{ALG@line}{step}{steps}
	\crefname{enumeratpropi}{property}{properties}
	\crefname{enumeratpropii}{property}{properties}
	\counterwithout{algorithm}{section}
	\renewcommand{\thealgorithm}{\arabic{algorithm}}

\newcounter{algsubstate}
\renewcommand{\thealgsubstate}{\alph{algsubstate}}

	\newcommand{\LL}{\mathcal L}

	\newcommand{\M}{M}
	\newcommand{\A}{\mathcal A}
	\newcommand{\B}{\mathcal B}

	\newcommand{\J}{{\rm J}}%

	\newcommand\sym[1]{\mathop{\mathbb{S}}^{#1}}

	\DeclareMathOperator{\blkdiag}{blkdiag}

	\DeclareMathOperator{\vecc}{vec}

	\newcommand{\bx}{\bm x}
	\newcommand{\bw}{\bm w}
	\newcommand{\bs}{\bm c}
	\newcommand{\bq}{\bm q}
	\newcommand{\bl}{\bm \lambda}

	\newcommand{\Ell}{j}

	\pgfplotsset{compat=1.16}

	\newcommand{\TheTitle}{Neural Network Training as an Optimal Control Problem\texorpdfstring{\\--- }{: }An Augmented Lagrangian Approach\texorpdfstring{ ---}{}}
	
	\newcommand{\TheAuthor}{%
		Brecht Evens,
        Puya Latafat,
		Andreas Themelis,
		Johan Suykens,
        and
		Panagiotis Patrinos%
	}%
	\newcommand{\TheFunding}{%
        This work was supported by the Research Foundation Flanders (FWO) research projects G0A0920N, G086518N, G086318N, and PhD grant 1196820N;
        Research Council KU Leuven C1 project No. C14/18/068;
        Fonds de la Recherche Scientifique -- FNRS and the Fonds Wetenschappelijk Onderzoek -- Vlaanderen under EOS project no 30468160 (SeLMA). Johan Suykens and Panagiotis Patrinos are affiliated to Leuven.AI - KU Leuven institute
		for AI, B-3000, Leuven, Belgium.
    }
	\newcommand{\TheKeywords}{%
		Neural networks,
		augmented Lagrangian method,
		Gauss-Newton method,
		dynamic programming%
	}
	
	\newcommand{\TheAbstract}{%
		Training of neural networks amounts to nonconvex optimization problems that are typically solved by using backpropagation and (variants of) stochastic gradient descent.
		In this work we propose an alternative approach by viewing the training task as a nonlinear optimal control problem.
		Under this lens, backpropagation amounts to the sequential approach (single shooting) to optimal control, where the states variables have been eliminated.
		It is well known that single shooting may lead to ill conditioning, and for this reason the simultaneous approach (multiple shooting) is typically preferred.
		Motivated by this hypothesis, an augmented Lagrangian algorithm is developed that only requires an approximate solution to the Lagrangian subproblems up to a user-defined accuracy.
		By applying this framework to the training of neural networks, it is shown that the inner Lagrangian subproblems are amenable to be solved using Gauss-Newton iterations.
		To fully exploit the structure of neural networks, the resulting linear least squares problems are addressed by employing an approach based on forward dynamic programming.
		Finally, the effectiveness of our method is showcased on regression datasets.
	}

\title{%
	\texorpdfstring{%
		\LARGE\bf\TheTitle\thanks{\TheFunding}%
	}{\TheTitle}%
}
\author{%
	\texorpdfstring{%
		Brecht Evens\affil{1}\thanks{\affil{1}%
			KU Leuven, Department of Electrical Engineering ESAT-STADIUS -- %
			Kasteelpark Arenberg 10, bus 2446, B-3001 Leuven, Belgium
			\newline
			{\sf
				\{%
					\href{mailto:brecht.evens@kuleuven.be}{brecht.evens},%
					\href{mailto:puya.latafat@kuleuven.be}{puya.latafat},%
					\href{mailto:johan.suykens@kuleuven.be}{johan.suykens},%
					\href{mailto:panos.patrinos@kuleuven.be}{panos.patrinos}%
				\}%
				\href{mailto:brecht.evens@kuleuven.be,puya.latafat@kuleuven.be,panos.patrinos@kuleuven.be}{@kuleuven.be}%
			}%
		}
        \and
        Puya Latafat\affil{1}
		\and
		Andreas Themelis\affil{2}\thanks{\affil{2}%
			Kyushu University, Faculty of Information Science and Electrical Engineering (ISEE) -- %
			744 Motooka, Nishi-ku 819-0395, Fukuoka, Japan
			\newline
			{\sf
				\href{mailto:andreas.themelis@ees.kyushu-u.ac.jp}{andreas.themelis@ees.kyushu-u.ac.jp}%
			}%
		}
		\and
		Johan Suykens\affil{1}
		\and
		Panagiotis Patrinos\affil{1}
	}{%
		\TheAuthor
	}%
}

\begin{document}

	\maketitle
	\begin{abstract}\TheAbstract \end{abstract}


	\section{Introduction}\label{sec:intro}
Feedforward deep neural networks (DNNs) are a prominent model for supervised learning, having a lot of success in various fields.
The primary objective of this work is to devise a novel method for training DNNs with smooth activation functions; this task can be formally stated as follows.

\begin{problem*}
	Given pairs
	\(
		\set*{
			(a^{(\ell)},b^{(\ell)})\in\R^{d_0}\times\R^{d_{N+1}}
		}_{\ell\in[m]}
	\),
	continuously differentiable functions
	\(
		\set*{
			\func{\Phi_j}{\R^{d_j}}{\R^{d_j}}
		}_{j\in[N+1]}
	\)
	(operating in an element-wise fashion), and $\mu_w>0$, find
	\(
		\set*{
			W_j\in\R^{d_j\times d_{j-1}}
		}_{j\in[N+1]}
	\)
	solutions to
	\begin{subequations}\label{subeq:NN}
		\begin{align}
		\nonumber
			\minimize_{W_1,\hdots,W_{N+1}}{}
		&
			\textstyle
			\tfrac{1}{2m}\sum_{\ell=1}^{m}\|\Phi_{N+1}(W_{N+1}x_N^{(\ell)})-b^{(\ell)}\|^2
			{}+{}
			\tfrac{\mu_w}{2}\sum_{j=1}^{N+1}\|W_j\|_F^2
		\\
		\label{eq:NN:dyn:init}
			\text{where}\;
		&
			x_0^{(\ell)}\coloneqq a^{(\ell)},\;
			\msmall{\ell\in[m]},
		\\
		\label{eq:NN:dyn}
		&
			x_j^{(\ell)}\coloneqq\Phi_j(W_jx_{j-1}^{(\ell)}),\;
			\msmall{j \in [N]},\;
			\msmall{\ell\in[m]}.
		\end{align}
	\end{subequations}
\end{problem*}

Here, \((a^{(\ell)},b^{(\ell)})\) are (given) \emph{training pairs}, \(N\in\N\) is the number of \emph{layers} of the network, each one having $d_i$ many \emph{neurons/nodes} and with $\Phi_i$ being the corresponding \emph{activation function}, and $\mu_w$ is a regularization parameter for the \emph{weights} \(W_i\) commonly used to avoid overfitting \cite{krogh1991simple}.
These optimization problems are typically solved using backpropagation \cite{rumelhart1986learning} along with (variants of) stochastic gradient descent, due to their simplicity and effectiveness.
However, these optimization methods suffer from various issues related to the challenging, highly nonconvex nature of the training task.
First and foremost, due to the prominence of local minima and saddle points, trained DNN models tend to generalize poorly to test data.
To alleviate this issue, various regularization methods have been introduced such as weight decay \cite{krogh1991simple}, batch normalization \cite{ioffe2015batch}, and dropout \cite{srivastava2014dropout}, typically reducing the overfitting of the training data.
More fundamentally, gradient-based methods are known to suffer from the vanishing gradient phenomenon \cite{hochreiter2001gradient}, where the gradients in the ouput layers of DNNs decrease exponentially with the number of layers.
Although recent studies have shown that piecewise affine activation functions such as ReLU, leaky ReLU \cite{maas2013rectifier}, and maxout unit \cite{goodfellow2013maxout} reduce the vanishing gradient problem by making the problem more sparse, the issue nevertheless persists especially in very deep networks.

To address these issues, in recent years a host of auxilary variable methods have been introduced where the network structure is represented by equality constraints and the space of learning parameters is extended.
By lifting the number of variables, these methods decompose the training task into a series of local subproblems which can be solved deterministically, typically using block coordinate descent (BCD) \cite{carreira2014distributed,zhang2017convergent, gu2020fenchel} or the alternating direction method of multipliers (ADMM) \cite{zhang2016efficient, taylor2016training, wang2019admm}.
BCD and ADMM have been successful for this task due to their ability to convert the equality constrained optimization problems into unconstrained problems, which can then be solved more efficiently than their constrained counterparts.
By increasing the dimension of the training problem, auxilary variable methods are able to alleviate some of the issues from which classical gradient-based methods suffer.
Most notably, it is observed that the vanishing gradient issue is alleviated as the auxilary variables circumvent long-term dependencies between the network weights during training \cite{zhang2016efficient}.
On the other hand, the increased dimensionality naturally makes the training task more challenging than when using classical gradient-based approaches.

The difference between traditional methods and auxilary variable methods can be related to concepts from optimal control by viewing the training task as a nonlinear optimal control problem.
Under this lens, auxilary vairable methods amount to the simultaneous approach (multiple shooting), whereas backpropagation amounts to the sequential approach (single shooting), where the state variables are eliminated \cite{lecun1988theoretical}.
As it is well known that single shooting may lead to ill conditioning of the optimization problem, it can be expected that multiple shooting methods can provide major advantages in the learning process of DNNs.

Motivated by this hypothesis, we develop a training methodology for neural networks based on an augmented Lagrangian framework that only requires finding approximate stationary points of the Lagrangian subproblems up to a user-defined accuracy.
To fully exploit the structure of feedforward neural networks, we additionally provide a computationally efficient approach to solve the inner subproblems based on forward dynamic programming.
The overall approach leads to an efficient and provably convergent methodology for solving the highly nonconvex optimization problems emerging in the neural network training task.

		\subsection{Contributions}
The contribution of this paper is twofold:
\begin{enumerate}[%
	leftmargin=0pt,
	labelwidth=7pt,
	itemindent=\labelwidth+\labelsep,
	label=\rlap{{\bf\arabic*)}}\hspace*{\labelwidth},
]
\item
	We introduce a novel augmented Lagrangian framework (ALM) for solving general nonconvex and nonsmooth equality constrained optimization problems.
	The framework is inspired by and extends \cite[Alg. 1]{grapiglia2020complexity} by waiving smoothness assumptions and relaxing the penalty update rule, yet preserving convergence to approximate KKT points in finite time.
\item
	We apply this framework to the training of DNNs, which we address from an optimal control perspective.
	The resulting optimization problem's structure has a twofold benefit: first, the inner Lagrangian subproblems are amenable to be addressed with fast methods such as Gauss-Newton (GN); in turn, forward dynamic programming (FDP) can conveniently be employed to efficiently solve the resulting linear least squares problems.
\end{enumerate}
To reflect the modularity and the contribution of each component, the three procedures (outer ALM, inner GN, and FDP) are outlined in three standalone algorithms, each addressing a dedicated general problem.

        \subsection{Organization}
The paper is organized as follows.
The notation is introduced in the next subsection.
An optimal control reformulation for the NN problem is presented in \cref{sec:OC}.
In \cref{sec:outerALM} a novel augmented Lagrangian method (ALM) is proposed for general equality constrained nonlinear programs.
The ALM method is specialized for training of neural networks with smooth activation functions in \cref{sec:innerGN}, where a procedure based on the Gauss-Newton method  and forward dynamic programming is proposed.
The proofs of all the results are deferred to the appendix.
Finally, numerical simulations showcasing the effectiveness of our proposed methodology on regression datasets are discussed in \cref{sec:experiments}.

		\subsection{Notation}
We use $[N]$ to denote the set of indices $\set{1,\dots,N}$.
We denote by $\R^n$ the standard $n$-dimensional Euclidean space with inner product $\innprod{{}\cdot{}}{{}\cdot{}}$ and induced norm $\|{}\cdot{}\|$.
The set of extended real numbers is defined as $\Rinf\coloneqq\R\cup\set{\infty}$, and we say that an extended-real valued function \(\func{f}{\R^n}{\Rinf}\) is proper if \(\dom f\coloneqq\set{x\in\R^n}[f(x)<\infty]\) is nonempty.
The set of real $n$-by-$m$ matrices is denoted by $\R^{n\times m}$.
Given \(A\in\R^{n\times m}\), $\|A\|_F$ is its Frobenius norm and \(\vecc(A)\in\R^{nm}\) is the vector obtained by stacking the columns of $A$ on top of one another.
The sets of symmetric, symmetric positive semi-definite and symmetric positive definite $n$-by-$n$ matrices are denoted by $\sym{n}$, $\sym{n}_+$ and $\sym{n}_{++}$, respectively.
For $V\in\sym{n}_{++}$ we define the scalar product $\innprod{x}{y}_V=\innprod{x}{Vy}$ and the induced norm $\|x\|_V=\sqrt{\innprod xx_V}$.
The \(n\)-by-\(n\) identity matrix is denoted by $\I_n$, or simply $\I$ when no ambiguity occurs.
The vector of all zeros with dimention $n$, and the $n$-by-$m$ matrix of all zeros are denoted by $0_n$, and $0_{n\times m}$, respectively.
The matrix Kronecker product is denoted by $\otimes$.
The Jacobian of a differentiable function \(\func{F}{\R^n}{\R^m}\), is denoted by \(\func{\J F}{\R^n}{\R^{m\times n}}\); \(\J_xF\) is a short-hand notation for the partial derivative \(\dep Fx\).

	\section{An optimal control reformulation}\label{sec:OC}
In the traditional approach, \eqref{eq:NN:dyn:init} and \eqref{eq:NN:dyn} are absorbed into the cost, thus forming an unconstrained minimization which is then solved by employing a stochastic (sub)gradient-type method.
Here we take an alternative approach by viewing the minimization as an optimal control problem with $N$ stages.
To this end, \eqref{subeq:NN} represents the dynamics of the problem and may compactly be written as
\begin{equation}\label{eq:NN:dyn:states}
    X_j = \Phi_j(W_j X_{j-1}),
\quad
	\msmall{j\in[N]}
\end{equation}
where $X_j \in \R^{d_j \times m}$ is a matrix whose $i$-th column is the vector $x_j^{(i)}$, for $i\in[m]$.
By similarly letting $A\in\R^{d_0\times m}$ and $Y\in\R^{d_{N+1}\times m}$ denote the  input and output matrices (constructed using vectors $a^{(i)}$, $b^{(i)}$), the following compact reformulation of \eqref{subeq:NN} is obtained
\begin{align*}
	\minimize_{\mathclap{\seq{W_i}[i\in{[N+1]}], \seq{X_i}[i\in{[N]}]}}{}
&\quad
	\textstyle
	\tfrac{1}{2m}\|\Phi_{N+1}(W_{N+1}X_N)-Y\|_F^2
	{}+{}
	\tfrac{\mu_w}{2}\sum_{i=1}^{N+1}\|W_i\|_F^2
\\
\numberthis\label{eq:OCP}
	\stt{}
&\quad
	\fillwidthof[r]{X_{j+1}}{X_0} = A
\\
&\quad
	X_{j+1}=\Phi_j(W_jX_{j-1}),\;
	\msmall{j\in[N]}.
\end{align*}

		\subsection{Vectorized form}\label{sec:vectorized}
For simplicity of exposition and computational convenience, we condense the optimization variables $W_i$ and $X_i$ into a single long vector \(\bm z=(\bm w, \bm x)\) with
\[
	\bm w
{}={}
	(\bm w_1,\dots,\bm w_{N+1})
~~\text{and}~~
	\bm x
{}={}
	(\bm x_1,\ldots, \bm x_N)
\]
where, letting $w_{i,j}\in\R^{d_{i-1}}$ denote the $j$-th row of $W_i$ and $x_i^{(j)}$ the $j$-th column of $X_i$ as in \cref{sec:intro},
\[
	\bm w_i
{}={}
	(w_{i,1},\dots,w_{i, d_i})\in\R^{d_id_{i-1}}
~\text{and}~
	\bm x_i
{}={}
	(x_i^{(1)}, \dots, x_i^{(m)})\in\R^{md_i}.
\]
In the vectorized notation, the cost function and the nonlinear constraints in \eqref{eq:OCP} may be represented by $f$ and $F(\bm z)=0$ with
\begin{align*}
	f(\bm z)
{}={} &
    \tfrac{1}{2m}
    \|H_{N+1}(\bw_{N+1},\bx_N) - \bm y\|^2
    {}+{}
    \tfrac{\mu_w}{2}\|\bw\|^2,
\\
    F(\bm z)
{}={} &
	\Big(
		\bm x_1-H_1(\bw_1, \bx_0), \ldots,
		\bm x_{N+1}-H_{N+1}(\bw_{N+1}, \bx_N)
	\Big),
\end{align*}
where \(\bm y=\vecc(Y)\) and
\begin{equation}\label{eq:H}
	H_j(\bw_j, \bx_{j-1})
{}={}
	\bigl(\Phi_j(W_jx_{j-1}^{(1)}),\dots,\Phi(W_jx_{j-1}^{(m)})\bigr).
\end{equation}

	\section{The outer ALM algorithm}\label{sec:outerALM}
With vectorized notation being adopted and as long as the activation functions \(\Phi_j\) are locally Lipschitz, the minimization in \eqref{eq:OCP} falls into the following general setting.

\begin{problem}[General ALM framework]\label{prob:ALM}%
	For a proper, lower semicontinuous, lower bounded \(\func{f}{\R^n}{\Rinf}\) and a locally Lipschitz \(\func{F}{\R^n}{\R^p}\) such that \(\set{\bm z\in\dom f}[F(\bm z)=0]\neq\emptyset\),
	\begin{equation}\label{eq:ALMP}
	\textstyle
		\minimize_{\bm z\in\R^n}f(\bm z)
	\quad\stt{}
		F(\bm z)=0.
	\end{equation}
\end{problem}

This section proposes a conceptual algorithm for addressing \cref{prob:ALM}, conceptual in the sense that, at this stage, no hint is given as to how the inner subproblems it involves can be solved.
The algorithm will be concretized in the subsequent \cref{sec:innerGN}, where an implementable procedure for addressing these inner steps is detailed.
The chosen method for the inner subproblems will ultimately require some additional structure and differentiability assumptions, which are nevertheless not needed for the (outer) ALM scheme presented in this section.
For the sake of generality of the discussion and to well pinpoint where each requirement is invoked, the convergence proof of the outer scheme is given in this broader setting.

Equality constrained minimization problems as \eqref{eq:ALMP} are amenable to be addressed by means of augmented Lagrangian methods.
For \(\beta>0\), we denote the corresponding \(\beta\)-augmented Lagrangian as
\begin{align}
	\LL_\beta(\bm z,\bm\lambda)
{}\coloneqq{} &
	f(\bm z)
	{}+{}
	\innprod{\bm\lambda}{F(\bm z)}
	{}+{}
	\tfrac\beta2\|F(\bm z)\|^2 \nonumber
\\
\label{eq:L'}
{}={} &
	f(\bm z)
	{}+{}
	\tfrac\beta2\|F(\bm z)+\bm\lambda\nicefrac{}{\beta}\|^2
	{}-{}
	\tfrac{1}{2\beta}\|\bm\lambda\|^2,
\end{align}
\begin{subequations}\label{eq:KKT}%
	and we say that \((\bm z,\bm\lambda)\) is an \(\varepsilon\)-KKT pair if
	\begin{align}
	\label{eq:KKTopt}
		\|\nabla_{\bm z}\LL(\bm z,\bm\lambda)\|_\infty
	{}\leq{} &
		\varepsilon,
	\quad\text{and}
	\\
	\label{eq:KKTfeas}
		\|F(\bm z)\|_\infty
	{}\leq{} &
		\varepsilon,
	\end{align}
	where \(\LL\coloneqq\LL_0\) is the (non-augmented) Lagrangian, and \(\nabla_{\bm z}\) denotes the gradient with respect to \(\bm z\) or, in case of lack of differentiability, any vector in the subdifferential \(\partial_{\bm z}\LL(\bm z,\bm\lambda)\).
\end{subequations}

Largely inspired by \cite[Alg. 1]{grapiglia2020complexity}, \cref{alg:ALM} hinges on the upper boundedness of the augmented Lagrangian along the iterates (see, e.g. \cite[Ex. 4.12]{birgin2014practical}) ensured by the initialization at a feasible point \(\bm z^0\).
Being not concerned with the tight rate analysis of \cite{grapiglia2020complexity}, we reduced the assumptions to the general setting of \cref{prob:ALM} and proposed a less conservative update rule for the penalty parameter.

\begin{algorithm}[t]
	\caption{ALM for \protect\cref{prob:ALM}}%
	\label{alg:ALM}%
\begin{algorithmic}[1]%
\setlength\itemsep{0.5ex}%
\Require
	\begin{tabular}[t]{@{}l@{}}%
		Initial feasible point \(\bm z^0\in\dom f\) s.t. \(F(\bm z^0)=0\),
	\\
		multiplier \(\bm\lambda^0\), and penalty \(\beta_0>0\)
	\\
		Parameters \(0<\gamma<1<\alpha,\xi\) and tolerance \(\varepsilon>0\)
	\end{tabular}
\item[For \(k=0,1,2\ldots\)]
\State\label{state:ALM:barx}%
	Set \(\hat{\bm z}^k=\bm z^k\) if \(\LL_{\beta_k}(\bm z^k,\bm\lambda^k)\leq\overbracket*{\smash{f(\bm z^0}{})}^{\mathclap{\LL_{\beta_k}(\bm z^0,\bm\lambda^k)}}\),
	or \(\hat{\bm z}^k=\bm z^0\) otherwise
\State\label{state:ALM:x}%
	\parbox[t]{0.98\linewidth}{%
		Starting at \(\hat z^k\), apply a descent method to compute an \(\varepsilon\)-stationary point \(\bm z^{k+1}\) of
		\begin{equation}\label{eq:ALM:min}
			\minimize{}\LL_{\beta_k}({}\cdot{},\bm\lambda^k),
		\end{equation}
		\ie a point \(\bm z^{k+1}\) such that
		\begin{equation}\label{prop:ALM:grad}
			\|\nabla_{\bm z}\LL_{\beta_k}(\bm z^{k+1},\bm\lambda^k)\|_\infty\leq\varepsilon
		\end{equation}
	}%
\State\label{state:ALM:y}%
	Set \(\bm\lambda^{k+1}=\bm\lambda^k+\beta_kF(\bm z^{k+1})\)
\If{ \(\|F(\bm z^{k+1})\|_\infty\leq\varepsilon\) }%
	\Statex\hspace*{\algorithmicindent}%
		{\algfont{Return}} \(\varepsilon\)-KKT pair \((\bm z^{k+1},\bm\lambda^{k+1})\)
\EndIf
\State\label{state:ALM:beta}%
	\begin{tabular}[t]{@{}r@{~~}l@{}}
		Set & \(\beta_{k+1}=\beta_k\) ~if \(\|F(\bm z^{k+1})\|_\infty\leq\gamma\|F(\bm z^k)\|_\infty\),\\
		or & \(\beta_{k+1}=\max\set{\xi\beta_k,\beta_0(k+1)^\alpha}\) ~otherwise.
	\end{tabular}
\end{algorithmic}
\end{algorithm}

\begin{thm}\label{thm:ALM}
	Applied to \cref{prob:ALM}, \cref{alg:ALM} terminates in finite time and yields an \(\varepsilon\)-KKT pair for \eqref{eq:ALMP}.
\end{thm}

Note that the result can cope with rather general functions \(f\) and \(F\), not necessarily derived from formulations as in \eqref{eq:OCP}.
The optimal control structure will instead be exploited in the following \cref{sec:innerGN} where an iterative method for addressing the inner problems at \cref{state:ALM:x} will be given.

\begin{rem}
	When \(f\) is lower bounded, then so is $\LL_\beta({}\cdot{},\bm\lambda)$ for any $\bm\lambda$, thus ensuring the existence of \(\varepsilon\)-stationary points as required in \cref{state:ALM:x} for any \(\varepsilon>0\).
	In the setting of the optimal control problem \eqref{eq:OCP}, not only is this condition trivially satisfied, but a feasible starting point \(\bm z^0\) can be obtained at virtually no cost by initializing the weights \(W_j^0\) and unrolling the dynamics to generate the \emph{state} variables \(X_j^0\).%
\end{rem}

	\section{The Lagrangian subproblem via Gauss-Newton iterations}\label{sec:innerGN}
In this section we present a procedure for solving the inner minimization \eqref{eq:ALM:min} in the setting of NNs with continuously differentiable activation functions.
With the notational conventions of \cref{sec:vectorized}, for a fixed multiplier \(\bl=(\bl_1,\dots,\bl_N)\) (\(\bl_j\) being the one associated with the \(j\)-th dynamics) the Lagrangian subproblem associated with \eqref{eq:OCP} is cast as follows.

\begin{problem}[Lagrangian subproblem]\label{prob:LagSub:Vec}%
	Given smooth functions \(\set*{H_j}_{j\in[N+1]}\), vectors \(\bl,\bm y\) of suitable sizes, and \(\beta,\mu_w>0\),%
	\begin{align*}
	\nonumber
		\minimize_{\bm z=(\bm w, \bm x)}{}
		\LL_\beta(\bm z,\bm\lambda)
	{}\coloneqq{}
	&
		\tfrac{1}{2m}\|H_{N+1}(\bm w_{N+1},\bm x_N)-\bm y\|^2
		{}+{}
		\tfrac{\mu_w}{2}\|\bm w\|^2
	\\
	&
		\textstyle
		\mathllap{
			-\tfrac{1}{2\beta}\|\bl\|^2
		}
		{}+{}
		\tfrac{\beta}{2}\sum_{j=1}^N\|\bm x_j- H_j(\bm w_j, \bm x_{j-1}) + \bl_j\nicefrac{}{\beta}\|^2.
	\end{align*}
\end{problem}

This smooth unconstrained least-squares problem is amenable to be solved by the Gauss-Newton (GN) method, which amounts to iteratively solving minimizations obtained after linearizing functions $H_j$ around the last iterates, and then applying a standard linesearch to guarantee convergence.
In the next subsection we derive explicit expressions of the Jacobian matrices involved in the linearization.

		\subsection{Gauss-Newton linearization and update direction}\label{sec:ABc}

Let  $\func{D^{\Ell}}{\R^{d_{{\Ell}-1}}}{\R^{d_{\Ell}\times d_{\Ell}}}$ be given by
\begin{align*}
	D^{\Ell}(v)
{}\coloneqq{} &
	\diag\bigl(
		\Phi_{\Ell}'(\innprod*{w_{{\Ell},1}}{v}),\dots,\Phi_{\Ell}'(\innprod*{w_{{\Ell},d_{\Ell}}}{v})
	\bigr),
\shortintertext{%
	(recall that \(\Phi_j\) operates element-wise) and define
}
	\mathcal D^{\Ell}
{}\coloneqq{} &
	\blkdiag\bigl(
		D^{\Ell}({x_{{\Ell}-1}^{(1)}}),\dots,D^{\Ell}({x_{{\Ell}-1}^{(m)}})
	\bigr).
\end{align*}
The Jacobians
\(
	\J_{\bm x_{{\Ell}-1}}H_j
{}\in{}
	\R^{md_{\Ell}\times md_{{\Ell}-1}}
\)
and
\(
	\J_{\bm w_{\Ell}}H_{\Ell}
{}\in{}
	\R^{md_{\Ell}\times d_{\Ell}d_{{\Ell}-1}}
\)
are then given by
\begin{align*}
	\J_{\bm x_{\Ell-1}}H_{\Ell}(\bm w_{{\Ell}}, \bm x_{\Ell-1})
{}={} &
	\mathcal D^{{\Ell}}\bigl(\I_m\otimes W_{{\Ell}}\bigr),
~~
	\text{and}
\\
	\J_{\bm w_{\Ell}}H_{\Ell}(\bm w_{\Ell}, \bm x_{{\Ell}-1})
{}={} &
	\mathcal D^{\Ell}\bigl(
		\I_{d_{\Ell}}\otimes x_{{\Ell}-1}^{(1)},\dots,\I_{d_{\Ell}}\otimes x_{{\Ell}-1}^{(m)}
	\bigr)^\top.
\end{align*}
If \(\bm z^l=(\bm w^l,\bm x^l)\) is the \(l\)-th iterate of a GN algorithm, denoting%
\begin{align*}
	A_{{\Ell}+1}
{}={} &
	\begin{ifcases}
		0_{md_1\times md_0} & \msmall{j=0}
	\\[3pt]
		\J_{\bm x_{\Ell}}H_{\Ell+1}(\bm w_{{\Ell}+1}^l,\bm x_{\Ell}^l) & \msmall{j\in[N]}
	\end{ifcases}
\\
\numberthis\label{eq:s}
	B_{\Ell}
{}={} &
	\J_{\bm w_{\Ell}}H_{\Ell}(\bm w_{\Ell}^l, \bm x_{{\Ell}-1}^l),
	\quad
	\msmall{\Ell\in[N+1]}
\\
	\bs_{\Ell}
{}={} &
	\begin{ifcases}
		H_{\Ell}(\bm w_{{\Ell}}^l,\bm x_{{\Ell}-1}^l)-A_{\Ell}\bm x_{{\Ell}-1}^l-B_{\Ell}\bm w_{\Ell}^l-\tfrac{1}{\beta}\bl_{\Ell}
		&
		\msmall{j\in[N]}
	\\[3pt]
		H_j(\bm w_j^l,\bm x_{j-1}^l)-A_j\bm x_{j-1}^l-B_j\bm w_j^l - \bm y
		&
		\msmall{j=N+1}
	\end{ifcases}
\end{align*}
the linearized minimization yielding the \(l\)-th GN update direction reduces to the following problem.

\begin{problem}[GN direction]\label{prob:GN}%
	Given \(\A = (A_1,\ldots, A_{N+1})\), \(\B = (B_1,\ldots, B_{N+1})\) and \(\bs=(\bs_1,\dots,\bs_{N+1})\) with matrices \(A_{\Ell},B_{\Ell}\) and vectors \(\bs_{\Ell}\) of suitable sizes, and given scalars $\beta,\mu_w>0$,
	\[
		\minimize_{\bm z=(\bw, \bx)}{
			\mathcal G_{\A, \B, \bs}(\bm z)
		}
		\quad
		\text{where,}
	\]
	denoting \((\delta_{\Ell},\rho_{\Ell})=(1,\beta)\) for \(\Ell\leq N\) and \((0,\frac1m)\) otherwise,%
	\[\mathtight[0.7]
		\textstyle
		\mathcal G_{\A, \B, \bs}(\bm z)
	{}\coloneqq{}
		\sum_{\Ell=1}^{N+1}\left(
			\tfrac{\rho_{\Ell}}{2}
			\|\delta_{\Ell}\bx_{\Ell} - A_{\Ell}\bx_{\Ell-1} - B_{\Ell}\bw_{\Ell} - \bs_{\Ell}\|^2
			{}+{}
			\tfrac{\mu_w}{2}
			\|\bw_{\Ell}\|^2
		\right).
	\]
\end{problem}

		\subsection{The Gauss-Newton algorithm}
The structure of \cref{prob:GN} emphasizes how variables are weakly coupled, a phenomenon that owes to the stagewise structure of the optimal control problem \eqref{eq:OCP}.
As a result, in spite of the large scale, \cref{prob:GN} admits a closed form solution that is efficiently retrievable with a forward dynamic programming (FDP) approach detailed in the following \cref{sec:FDP}.
This routine may then be invoked by the GN method, synopsized in \cref{alg:GN}, when computing the update directions at \cref{state:GN:min}.

\begin{algorithm}[t]
	\caption{Gauss-Newton procedure for \protect\cref{prob:LagSub:Vec}}%
	\label{alg:GN}%
\begin{algorithmic}[1]%
\setlength\itemsep{0.5ex}%
\Require
	Initial point
	\(
		\bm z^0
	{}={}
		(\bm w^0, \bm x^0)
	\)
	and \(0<\eta_1,\eta_2<1\)
\item[For \(l=0,1,2\ldots\)]%
\State\label{state:GN:min}%
	\begin{tabular}[t]{@{}l@{}}
		{\sc [update direction]}~~
		set \({\bm p}^l=\bar{\bm z}^l-\bm z^l\), where \(\bar{\bm z}^l=(\bar{\bm w}^l,\bar{\bm x}^l)\)
	\\
		solves \cref{prob:GN} with $\A,\B,\bs$ as in \eqref{eq:s}
	\end{tabular}
\State\label{state:GN:LS}%
	\begin{tabular}[t]{@{}l@{}}
		{\sc \fillwidthof[l]{[update direction]}{[linesearch]}}~~
		set \(\bm z^{l+1}=\bm z^l+\tau_l\bm p^l\), where \(\tau_l\)%
	\\
		is the largest number in \(\set*{1,\eta_1,\eta_1^2,\ldots}\)
		 such that
	\\[4pt]
		\quad
		\(\displaystyle
			\LL_\beta(\bm z^l + \tau_l\bm p^l,\bl)
		{}\leq{}
			\LL_\beta(\bm z^l,\bl)
			{}-{}
			\eta_2\tau_l \mathcal G_{\A,\B,0}(\bm p^l)
		\)%
	\\[4pt]
		with  \(\LL_\beta\) and \(\mathcal G_{\A,\B,\bs}\) as in \cref{prob:GN,prob:LagSub:Vec}
	\end{tabular}
\If{ \(\|\nabla_{\bm z}\LL_{\beta}(\bm z^{l+1},\bm\lambda)\|_\infty\leq\varepsilon\) }%
	\Statex\hspace*{\algorithmicindent}%
		{\algfont{Return}}  \(\bm z^{l+1}= (\bw^{l+1}, \bx^{l+1})\)
\EndIf
\end{algorithmic}
\end{algorithm}

In the next lemma we show that the GN method yields an $\varepsilon$-stationary solution for the original Lagrangian subproblem.

\begin{lem}\label{thm:GN:ST}
	Applied to \cref{prob:LagSub:Vec}, \cref{alg:GN} terminates in finite time yielding an $\varepsilon$-stationary solution.
\end{lem}

		\subsection{Forward dynamic programming}\label{sec:FDP}
In this subsection we propose a recursive procedure for solving \cref{prob:GN} with given matrices $A_{\Ell}\in\R^{r_{{\Ell}}\times r_{{\Ell}-1}}$, $B_{\Ell}\in\R^{r_{\Ell} \times s_{\Ell}}$, and vectors $\bs_{\Ell}\in\R^{r_{\Ell}}$, $\Ell\in[N+1]$, thus providing an efficient routine for \cref{state:GN:min} of \cref{alg:GN}.
Inspired by the idea of forward dynamic programming, the minimization may be split into a series of simpler subproblems that are solved in a recursive manner:
{\mathtight[0.8]%
	\makeatletter
		\let\oldEll\Ell
		\renewcommand{\Ell}{\@ifstar\@@Ell\@Ell}%
		\newcommand{\@Ell}{\oldEll-1}%
		\newcommand{\@@Ell}{\oldEll}%
	\makeatother
	\begin{align}
	\label{eq:FDR:V1}
		V_1^\star(\bx_1)
	{}={} &
		\min_{\bw_1}~\Bigl\{
			\tfrac{\rho_1}{2}
			\|\bx_1-A_1\bx_0-B_1\bw_1-\bs_1\|^2
			{}+{}
			\tfrac{\mu_w}2
			\|\bw_1\|^2
		\Bigr\}
	\\
	\nonumber
		V_{\Ell*}^\star(\bx_{\Ell*})
	{}={} &
		\min_{\mathclap{\bx_{\Ell},\bw_{\Ell*}}}~\Bigl\{
			V_{\Ell}^\star(\bx_{\Ell})
			{}+{}
			\tfrac{\rho_{\Ell*}}{2}
			\|\bx_{\Ell*}-A_{\Ell*}\bx_{\Ell}-B_{\Ell*}\bw_{\Ell*}-\bs_{\Ell*}\|^2
	\\
	\label{eq:FDR:Vk}
	&
			\hphantom{\min~\Bigl\{}
			{}+{}
			\tfrac{\mu_w}2
			\|\bw_{\Ell*}\|^2
		\Bigr\},
		\quad
		\msmall{\Ell*=2,\dots,N}
	\\
	\nonumber
		V_{N+1}^\star
	{}={} &
		\min_{\mathclap{\bx_N,\bw_{N+1}}}~\Bigl\{
			V_N^\star(\bx_N)
			{}+{}
			\tfrac{\rho_{N+1}}{2}
			\|A_{N+1}\bx_N + B_{N+1}\bw_{N+1} + \bs_{N+1}\|^2
	\\
	\label{eq:FDR:Vfinal}
	&
			\hphantom{\min~\Bigl\{}
			{}+{}
			\tfrac{\mu_w}2 \|\bw_{N+1}\|^2
		\Bigr\}.
	\end{align}%
}%
Each stage consists of minimization of the sum of  the cost at the current stage and the optimal cost from the previous stage.
The cost at the final stage $V_{N+1}^\star$ is equal to the optimal cost for \cref{prob:GN}.
In order to obtain closed form solutions for each of the above minimizations, let $E_{\Ell}\in\R^{s_{\Ell}\times r_{\Ell}}$ and $G_{\Ell},M_{\Ell},S_{\Ell}\in\R^{r_{\Ell}\times r_{\Ell}}$, $\Ell\in[N+1]$, be defined as
\begin{align}
\label{eq:E_G}
	E_{\Ell}
{}={} &
	\left(
		\tfrac{\mu_w}{\rho_{\Ell}}\I+B_{\Ell}^{\top}B_{\Ell}
	\right)^{-1}B_{\Ell}^\top,
\\
	G_{\Ell}
{}={} &
	\I - B_{\Ell} E_{\Ell},
\\
\label{eq:MatrixUpdates_M}
	M_{\Ell}
{}={} &
	\begin{ifcases}
		\frac1{\rho_1}\I+\frac{1}{\mu_w}B_1B_1^{\top}
	&
		\msmall{j=1}
	\\[4pt]
		\tfrac1{\rho_{\Ell}}\I
		{}+{}
		\tfrac{1}{\mu_w}B_{\Ell}B_{\Ell}^{\top}
		{}+{}
		A_{\Ell}
		M_{\Ell-1}A_{\Ell}^{\top}
	&
		\msmall{j>1}
	\end{ifcases}
\\
\label{eq:MatrixUpdates_S}
	S_j
{}={} &
	\begin{ifcases}
		\I_{r_0}
	&
		\msmall{j=1}
	\\
		M_{\Ell-1} - M_{\Ell-1} A_{\Ell}^\top M_{\Ell}^{-1}A_{\Ell} M_j.
		\otherwise[\msmall{j>1}.]
	\end{ifcases}
\end{align}
Note that matrices $S_{\Ell}$ need not be computed explicitly.
Instead, given a vector $v\in\R^{r_{\Ell}}$, $S_{\Ell}v$ is computed as follows:
\begin{align}\label{eq:S_kvUpdate}
	\begin{cases}[@{}c@{~}l@{}]
		\text{\it(i)} & \text{solve the linear system } M_{\Ell} \bar{v} = A_{\Ell}(M_{\Ell-1}{v})
	\\
		\text{\it(ii)} & \text{set } S_{\Ell}v=M_{\Ell-1} \big( v - A_{\Ell}^{\top} \bar{v} \big).
	\end{cases}
\end{align}

\begin{algorithm}[t]
	\caption{Recursive solution to \protect\cref{prob:GN} with FDP}%
	\label{alg:FDP}%
\let\OldComment\Comment
\renewcommand{\Comment}[1]{{\color{gray}\footnotesize\OldComment{#1}}}%

\begin{algorithmic}[1]
\setlength\itemsep{0.5ex}%
\Require
	\begin{tabular}[t]{@{}l@{}}%
		Initial state \(\bx_0\in\R^{r_0}\)
	\\
		set $M_1=\tfrac1{\rho_1}\I+\tfrac{1}{\mu_{u}}B_1B_1^{\top}$, $S_1=\I_{r_0}$, $\bq_0=\bx_0$
	\end{tabular}
\State\label{state:Forward}%
	{\sc [Forward recursion]}\quad
	For $\Ell=1,\ldots,N$%
	\begin{enumerate}[%
		leftmargin=3em,
		label=(\alph*),
		ref={step \protect{\ref{state:Forward}}(\alph*)},
	]
	\item\label{state:Forward:a}
		solve the linear system $M_{\Ell}\tilde{\bs}_{\Ell} = \bs_{\Ell}$,
	\item\label{state:Forward:b}
		$\tilde{\bq}_{\Ell}=S_{\Ell}\bq_{\Ell-1}$
	\Comment{as described in \eqref{eq:S_kvUpdate}}
	\item
		\(
			\bq_{\Ell} = \rho_{\Ell}G_{\Ell}A_{\Ell} \tilde{\bq}_{\Ell} + \tilde{\bs}_{\Ell}
		\)
	\item  
		\(
			M_{\Ell+1}
		{}={}
			\tfrac1{\rho_{\Ell+1}}\I+\tfrac{1}{\mu_{u}}B_{\Ell+1}B_{\Ell+1}^{\top}+A_{\Ell+1}M_{\Ell}A_{\Ell+1}^{\top}
		\)
	\end{enumerate}
\State\label{state:Backward}%
	{\sc [Backward recursion]}
		\begin{enumerate}[{label=(\alph*)},{ref={step \protect{\ref{state:Backward}}(\alph*)}}]
		\item\label{state:Backward:a}
			$\tilde{\bx}_{N+1} = - S_{N+1}\big(A_{N+1}^{\top}G_{N+1}\bm{c}_{N+1}\big)$
		\item[]
			$\tilde{\bq}_{N+1}=S_{N+1}\bq_{N}$
			\Comment{as described in \eqref{eq:S_kvUpdate}}
		\item
			\(
				\fillwidthof[r]{\bw_{N+1}}{\bx_N}
			{}={}
				\tilde{\bq}_{N+1} + \rho_{N+1}\tilde{\bx}_{N+1}
			\)
		\item[]
			\(
				\bw_{N+1}
			{}={}
				-E_{N+1}\big(A_{N+1}\bx_{N}+\bs_{N+1}\big)
			\)
		\item\label{state:Backward:c}
			For $\Ell=N,\ldots,2$:
			\begin{itemize}
			\item[]
				$\tilde{\bx}_{\Ell} = S_{\Ell}\big(A_{\Ell}^{\top}G_{\Ell}(\bx_{\Ell} - \bm{c}_{\Ell})\big)$
			\Comment{as described in \eqref{eq:S_kvUpdate}}
			\item[]
				\(
					\bx_{\Ell-1}
				{}={}
					\tilde{\bq}_{\Ell}+ \rho_{\Ell}\tilde{\bx}_{\Ell}
				\)
			\item[]
				\(
					\fillwidthof[r]{\bx_{\Ell-1}}{\bw_{\Ell}}
				{}={}
					E_{\Ell}\big(\bx_{\Ell}-A_{\Ell}\bx_{\Ell-1}-\bs_{\Ell}\big)
				\)
			\end{itemize}
		\item
			\(
				u_1
				{}={}
				E_1\big(\bx_1-A_1\bx_{0}-\bs_1\big)
			\)
		\end{enumerate}
\item[\algfont{Return}]
	$\bm z=(\bw, \bx)$ with $\bw=(\bw_1,\ldots,\bw_{N+1})$, $\bx=(\bx_1,\ldots,\bx_N)$
\end{algorithmic}
\end{algorithm}

The FDP procedure is presented in \cref{alg:FDP}.
Other than matrix-vector products, the algorithm requires solving linear systems several times, which may be performed by computing the \emph{Cholesky factorization} of $M_{\Ell}$ and $\tfrac{\mu_w}{\rho_{\Ell}}\I+B_{\Ell}^{\top}B_{\Ell}$ once, thus resulting in operations involving simple forward and backward substitution steps that substantially reduce the computational overhead.

\begin{rem}[Positive definiteness] \label{rem:RkUk}
	Since \(\rho_j,\mu_w>0\), \(G_{\Ell},M_{\Ell}\in\sym{r_{\Ell}}_{++}\) for any \(\Ell\).
	Furthermore, using the Woodbury matrix identity and \eqref{eq:MatrixUpdates_M}, the following alternative expression for $S_{\Ell+1}$ is obtained
	\begin{equation}\label{eq:SnM}
		\textstyle
		S_{\Ell+1}
	{}={}
		\bigl(M_{\Ell}^{-1} + \rho_{\Ell+1}A_{\Ell+1}^{\top}G_{\Ell+1}A_{\Ell+1}\bigr)^{-1},
	\end{equation}
	establishing that also $S_{\Ell+1}\in\sym{r_{\Ell}}_{++}$.
\end{rem}

The optimality of the solution obtained by the FDP procedure is established in the next lemma.

\begin{lem}\label{thm:FDP}
	 Suppose that $\mu_w>0$.
	 Then, $\bm z=(\bw, \bx)$ generated by \cref{alg:FDP} is the unique minimizer of \cref{prob:GN}.
\end{lem}

	\section{Numerical experiments}\label{sec:experiments}
\begin{table*}
	\caption{%
		Numerical results for training a three-layer network with varying input dimension and noise level using ALM, Adam and SGD.%
	}%
	\label{tab:benchmark}%
	\centering
	\begin{filecontents*}{Experiments/mean_table_three_layer_lowtolerance.csv}
S,d,epsilon,Training error,Test error,Feasibility error,Norm gradient,Lag evals,Jac evals,ALM iters,GN iters,Time (mm:ss),N
250,5,10\%,5.37e-2,5.04e-2,4.13e-4,2.58e-3,22,19,6,13,0:09,2.0
250,5,20\%,6.93e-2,6.62e-2,3.19e-4,3.10e-3,22,19,6,13,0:09,2.0
250,10,10\%,6.52e-2,6.47e-2,3.61e-4,3.60e-3,32,25,6,19,0:14,2.0
250,10,20\%,7.95e-2,8.08e-2,2.86e-4,3.91e-3,35,27,6,20,0:15,2.0
250,15,10\%,7.36e-2,7.98e-2,3.44e-4,4.99e-3,40,29,6,22,0:17,2.0
250,15,20\%,8.76e-2,9.49e-2,4.06e-4,4.57e-3,41,30,6,23,0:17,2.0
\end{filecontents*}
\pgfplotstabletypeset[
	string type,
	col sep=comma,
	columns={[index]1,[index]2,[index]3,[index]4,[index]7,[index]8,[index]9,[index]10,[index]11},
	display columns/0/.style={string type,column name={$d_0$}},
	display columns/1/.style={string type,column name={$\delta_0$}},
	display columns/2/.style={string type,column name={Training}},
	display columns/3/.style={string type,column name={Test}},
	display columns/4/.style={string type,column name={$\LL_{\beta_k}$}},
	display columns/5/.style={string type,column name={$\nabla_x \LL_{\beta_k}$}},
	display columns/6/.style={string type,column name={ALM}},
	display columns/7/.style={string type,column name={GN}},
	display columns/8/.style={string type,column name={Time}, string replace={0:}{b}},
	every head row/.style={
		before row={%
			\toprule
			&& \multicolumn{7}{c}{\textbf{ALM}}\\
			\cmidrule(lr){3-9}
		},
		after row/.add={}{
			&& MSE & MSE & evals & evals & iters & iters & (m:ss)\\ \midrule
			}
	},
	every last row/.style={after row=\bottomrule}
]{Experiments/mean_table_three_layer_lowtolerance.csv}
	\hspace{0.1cm}%
	\begin{filecontents*}{Experiments/mean_table_three_layer.csv}
S,d,epsilon,Training error,Test error,Feasibility error,Norm gradient,Lag evals,Jac evals,ALM iters,GN iters,Time (mm:ss),N,Adam training loss,Adam test loss,SGD training loss,SGD test loss,Time Adam,Time SGD
250,5,10\%,5.38e-2,5.05e-2,2.13e-4,4.61e-4,24,21,6,14,0:11,2.0,5.36e-2,5.05e-2,5.47e-2,5.12e-2,0:14,0:13
250,5,20\%,6.94e-2,6.62e-2,1.91e-4,5.94e-4,24,21,6,14,0:10,2.0,6.92e-2,6.63e-2,7.03e-2,6.68e-2,0:15,0:12
250,10,10\%,6.50e-2,6.44e-2,1.76e-4,4.64e-4,41,32,7,25,0:18,2.0,6.50e-2,6.44e-2,6.56e-2,6.48e-2,0:15,0:12
250,10,20\%,7.95e-2,8.08e-2,2.13e-4,5.53e-4,40,31,7,24,0:17,2.0,7.93e-2,8.07e-2,8.04e-2,8.15e-2,0:14,0:13
250,15,10\%,7.36e-2,7.98e-2,2.02e-4,5.53e-4,44,31,6,24,0:19,2.0,7.34e-2,7.96e-2,7.63e-2,8.35e-2,0:15,0:12
250,15,20\%,8.76e-2,9.48e-2,2.09e-4,6.26e-4,46,33,7,26,0:20,2.0,8.74e-2,9.47e-2,9.07e-2,9.93e-2,0:15,0:12
\end{filecontents*}
\pgfplotstabletypeset[
	string type,
	col sep=comma,
	columns={[index]13,[index]14,[index]17,[index]15,[index]16,[index]18},
	display columns/0/.style={string type,column name={Training}},
	display columns/1/.style={string type,column name={Test}},
	display columns/2/.style={string type,column name={Time}},
	display columns/3/.style={string type,column name={Training}},
	display columns/4/.style={string type,column name={Test}},
	display columns/5/.style={string type,column name={Time}},
	every head row/.style={
		before row={%
			\toprule
			\multicolumn{3}{c}{\textbf{Adam}} & \multicolumn{3}{c}{\textbf{SGD}}\\
			\cmidrule(lr){1-3}\cmidrule(lr){4-6}
		},
		after row/.add={}{
			MSE & MSE & (m:ss) & MSE & MSE & (m:ss)\\ \midrule
			}
	},
	every last row/.style={after row=\bottomrule}
]{Experiments/mean_table_three_layer.csv}
\end{table*}

		\subsection{Design of numerical experiments}
			We will generate training (and test) pairs $\set*{(a^{(\ell)},b^{(\ell)})}_{\ell\in[m]}$ for a three-layer neural network under the regression setting, analogous to the approach in \cite{cui2020multicomposite}, as follows:
\[
	b^{(\ell)} = W_3 \Phi(W_2 \Phi(W_1 a^{(\ell)})) + \delta
\]
where \(a^{(\ell)} \sim \mathcal{N}(\mu, \Sigma)\) and \(\delta \sim \delta_0 \mathcal{N}(0, 1)\).
The mean \(\mu \in \R^{d_0}\) and an additional random matrix \(\Sigma_0 \in \R^{d_0 \times d_0}\) are generated by a normal distribution with standard deviation 0.2, and the covariance \(\Sigma\) is set to be \(\Sigma_0^\top \Sigma_0\).
The three-layer network consists of $N=2$ hidden layers with respectively 20 and 5 neurons.
As activation function the softplus function is used, i.e. \(\Phi(x)\coloneqq\ln(1+\exp(x))\), a smooth approximation to the ReLU activation function which is often used in deep learning and known for its faster convergence.
The weights \(W_i\) of the neural network are initialized according to Kaiming \cite{he2015deep}, which is a weight initialization procedure suitable for networks consisting of softplus activation functions, and we obtain a feasible starting point $\bm z^0$ by applying \eqref{eq:NN:dyn:states} recursively.
All networks in this section are trained with regularization parameter $\mu_w = 0.1$.
The following parameters for \cref{alg:ALM} are used:
\begin{align*}
	\bm\lambda^0
{}={} &
	0,
&&
	\beta_0 = 0.001 f(\bm z^0),
&&
	\gamma = 0.5,
\\
    \alpha
{}={} &
	2,
&&
	\varepsilon = 10^{-3},
&&
	\xi = 2.
\end{align*}
Furthermore, to prevent solving the inner problems \eqref{eq:ALM:min} up to an unnecessarily high tolerance $\varepsilon$ in the first iterations, \cref{prop:ALM:grad} is relaxed as follows:
\begin{equation}\label{eq:updatingrule:epsilon}
	\|\nabla_{\bm z}\LL_{\beta_k}(\bm z^{k+1},\bm\lambda^k)\|_\infty
{}\leq{}
	\varepsilon_k
{}\coloneqq{}
    \max\left(\bar{\varepsilon},{}0.5\varepsilon_{k-1}\right)
\end{equation}
with \(\varepsilon_0 = 10^{-1}\) and $\bar{\varepsilon} = 10^{-2}$.
Finally, the following parameters for the line search in \cref{alg:GN} are used:
\begin{align*}
    \eta_1 &{}= 0.8, &&\eta_2 = 0.1.
\end{align*}
The ALM framework and corresponding Gauss-Newton procedure are implemented using the SciPy sparse matrix library \cite{virtanen2020scipy} in Python.
The CHOLMOD library \cite{chen2008algorithm} is used to factorize \((\frac{\mu_w}{\rho_j}I + B_j^\top B_j)\) and \(M_j\), which prevents the costly explicit computation of \(E_j\) and \(M_{j}^{-1}\).
All experiments are conducted on a HP elitebook 845 G7 with a 1.7GHz AMD Ryzen 7 PRO 4750U processor and 32 GB RAM.

		\subsection{Numerical results and discussion}
The left-hand side of \cref{tab:benchmark} shows the numerical results for training the previously introduced feedforward neural networks with varying input dimension $d_0$ and noise level $\delta_0$ (averaged over 15 simulations) using our proposed ALM method, which in a couple of ALM iterations yields an $\bar\varepsilon$-KKT pair (as \eqref{eq:KKTopt} is satisfied for $\bar\varepsilon$ instead of $\varepsilon$).

All experiments are performed with a fixed sample size $m=250$ for the training and test datasets.
We should remark that the current implementation does not scale well with the sample size $m$ both in terms of memory usage and computation time, as the matrices $M_j$ in the FDP procedure become increasingly large.
For this reason, our method would greatly benefit from a mini-batch implementation where the training set is split into smaller batches to compute the inner GN steps.
This is considered for future work.

\begin{figure}[H]
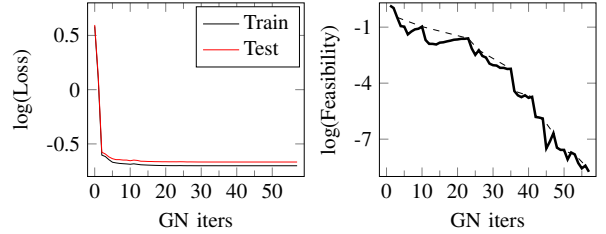

    \centering
    \includetikz{Loss}%
    \includetikz{Feas}%
    \caption{%
		Typical performance of the ALM algorithm.
		(Left) Training and test loss.
		(Right) Feasibility.
		The dashed line connects the points of the ALM iterations.%
	}%
	\label{fig:Typical performance ALM algorithm}%
\end{figure}

The typical performance of the ALM algorithm is visualized in \cref{fig:Typical performance ALM algorithm} for a simulation with $d_0 = 15$ and $\delta_0 = 20\%$ and tolerance $\varepsilon = 10^{-7}$ instead of $\varepsilon = 10^{-3}$.
In the earlier GN iterations mainly the loss is reduced, while in the final iterations the feasibility is recovered as the penalty parameter increases in the outer ALM iterations.
For this reason, it makes sense to terminate our algorithm at tolerance $\varepsilon = 10^{-3}$, as in neural network training we are mainly interested in reducing the loss.

		\subsection{Comparison with first-order methods}
We compare our previously obtained results with two commonly used first-order methods for stochastic optimization, namely Adam \cite{kingma2014adam} and stochastic gradient descent (SGD).
We use the default implementations of these algorithms provided by the Keras library using the TensorFlow \cite{abadi2016tensorflow} backend with batch size 10, MSE loss function and additional $\ell_2$ regularization with parameter $\mu_w$.

The right portion of \cref{tab:benchmark} shows the numerical results for training the three-layer network using Adam and SGD (averaged over 15 simulations) for 1000 epochs.
No early stopping or other monitoring callbacks are used, minimizing the computation time per epoch.
SGD is typically susceptible to stagnate at suboptimal points where it ceases to make significant progress, which explains its higher training MSE compared to Adam.
When comparing with Adam and SGD it can be seen that our method tends to converge towards very good local optima, surpassing the performance of SGD and occasionally even finding a better local minimum than Adam.
Furthermore, the computation time of our methodology for training the introduced networks is reasonably similar the ones of Adam and SGD.
Overall, these results are encouraging as our method is expected to greatly benefit from a mini-batch implementation, further reducing the computation time and increasing scalability.

    \section{Conclusions}\label{sec:conc}
In this paper a novel procedure for training of neural networks was introduced that leverages an optimal control view, and relies on three main components.
First, a novel augmented Lagrangian method is presented for general nonsmooth nonconvex equality constrained problems, which attains an $\varepsilon$-KKT solution in finite time.
Second, when applied to the DNN problem we propose to solve the Lagrangian subproblems by employing Gauss-Newton iterations resulting in a series of linear least squares problems.
Third, owing to the stagewise structure in the optimal control formulation, we solve the linear least squares GN problems through a simple recursive procedure based on forward dynamic programming.
We observed encouraging results in comparison to fast first-order solvers such as Adam which are often used in a heuristic manner without theoretical guarantees.
In the current implementation our method is not competitive when using large numbers of training data.
Future research directions include extending our scheme to mini-batch settings to tackle this issue.
It is also interesting to extend the framework to allow for nonsmooth activations functions.


	\bibliographystyle{plain}
	\bibliography{TeX/Bibliography.bib}

\begin{thebibliography}{10}

\bibitem{abadi2016tensorflow}
M.~Abadi, A.~Agarwal, et~al.
\newblock Tensorflow: Large-scale machine learning on heterogeneous distributed
  systems.
\newblock {\em arXiv:1603.04467}, 2016.

\bibitem{bertsekas2016nonlinear}
D.~Bertsekas.
\newblock {\em Nonlinear Programming}.
\newblock Athena Scientific, 2016.

\bibitem{birgin2014practical}
E.~Birgin and J.~Mart{'i}nez.
\newblock {\em Practical Augmented {L}agrangian Methods for Constrained
  Optimization}.
\newblock SIAM, 2014.

\bibitem{carreira2014distributed}
M.~Carreira-Perpinan and W.~Wang.
\newblock Distributed optimization of deeply nested systems.
\newblock In {\em Artif. Intell. Stat.}, pages 10--19, 2014.

\bibitem{chen2008algorithm}
Y.~Chen, T.~Davis, et~al.
\newblock Algorithm 887: {CHOLMOD}, supernodal sparse {Cholesky} factorization
  and update/downdate.
\newblock {\em ACM Trans Math Softw}, 35(3), oct 2008.

\bibitem{cui2020multicomposite}
Y.~Cui, Z.~He, and J.~Pang.
\newblock Multicomposite nonconvex optimization for training deep neural
  networks.
\newblock {\em SIAM J. Optim.}, 30(2):1693--1723, 2020.

\bibitem{facchinei2003finite}
F.~Facchinei and J.~Pang.
\newblock {\em Finite-dimensional variational inequalities and complementarity
  problems}, volume~II.
\newblock Springer, 2003.

\bibitem{goodfellow2013maxout}
I.~Goodfellow, D.~Warde-Farley, et~al.
\newblock Maxout networks.
\newblock In {\em Int. Conf. Mach. Learn.}, pages 1319--1327, 2013.

\bibitem{grapiglia2020complexity}
G.~Grapiglia and Y.~Yuan.
\newblock On the complexity of an augmented {L}agrangian method for nonconvex
  optimization.
\newblock {\em IMA J. Numer. Anal.}, jul 2020.

\bibitem{gu2020fenchel}
F.~Gu, A.~Askari, and L.~El~Ghaoui.
\newblock Fenchel lifted networks: A {L}agrange relaxation of neural network
  training.
\newblock In {\em Int. Conf. Artif. Intell. Stat.}, pages 3362--3371, 2020.

\bibitem{he2015deep}
K.~He, X.~Zhang, S.~Ren, and J.~Sun.
\newblock Deep residual learning for image recognition.
\newblock {\em arXiv:1512.03385}, 2015.

\bibitem{hochreiter2001gradient}
S.~Hochreiter, Y.~Bengio, et~al.
\newblock Gradient flow in recurrent nets: the difficulty of learning long-term
  dependencies, 2001.

\bibitem{ioffe2015batch}
S.~Ioffe and C.~Szegedy.
\newblock Batch normalization: Accelerating deep network training by reducing
  internal covariate shift.
\newblock {\em arXiv:1502.03167}, 2015.

\bibitem{kingma2014adam}
D.~Kingma and J.~Ba.
\newblock Adam: A method for stochastic optimization.
\newblock {\em arXiv:1412.6980}, 2014.

\bibitem{krogh1991simple}
A.~Krogh and J.~Hertz.
\newblock A simple weight decay can improve generalization.
\newblock In {\em Proc. 4th Int. Conf. NIPS}, page 950–957. Morgan Kaufmann
  Publishers Inc., 1991.

\bibitem{lecun1988theoretical}
Y.~LeCun.
\newblock {\em A theoretical framework for back-propagation}.
\newblock IEEE Computer Society Press, 1992.

\bibitem{maas2013rectifier}
A.~Maas, A.~Hannun, and A.~Ng.
\newblock Rectifier nonlinearities improve neural network acoustic models.
\newblock In {\em Proc ICML}, volume~30, page~3, 2013.

\bibitem{rockafellar2009variational}
R.~Rockafellar and R.~Wets.
\newblock {\em Variational Analysis}, volume 317.
\newblock Springer, 2009.

\bibitem{rumelhart1986learning}
D.~Rumelhart, G.~Hinton, and R.~Williams.
\newblock Learning representations by back-propagating errors.
\newblock {\em Nature}, 323(6088):533--536, 1986.

\bibitem{srivastava2014dropout}
N.~Srivastava, G.~Hinton, et~al.
\newblock Dropout: a simple way to prevent neural networks from overfitting.
\newblock {\em JMLR}, 15(1):1929--1958, 2014.

\bibitem{taylor2016training}
G.~Taylor, R.~Burmeister, et~al.
\newblock Training neural networks without gradients: A scalable {ADMM}
  approach.
\newblock In {\em Int. Conf. Mach. Learn.}, pages 2722--2731, 2016.

\bibitem{virtanen2020scipy}
P.~Virtanen, R.~Gommers, et~al.
\newblock {SciPy} 1.0: Fundamental algorithms for scientific computing in
  {Python}.
\newblock {\em Nature Methods}, 17:261--272, 2020.

\bibitem{wang2019admm}
J.~Wang, F.~Yu, et~al.
\newblock {ADMM} for efficient deep learning with global convergence.
\newblock In {\em Proc. 25th ACM SIGKDD Int. Conf. Knowledge Discov. Data
  Min.}, pages 111--119, 2019.

\bibitem{zhang2017convergent}
Z.~Zhang and M.~Brand.
\newblock Convergent block coordinate descent for training {Tikhonov}
  regularized deep neural networks.
\newblock {\em arXiv:1711.07354}, 2017.

\bibitem{zhang2016efficient}
Z.~Zhang, Y.~Chen, and V.~Saligrama.
\newblock Efficient training of very deep neural networks for supervised
  hashing.
\newblock In {\em Proc. IEEE Conf. Comput. Vision Pattern Recogn.}, pages
  1487--1495, 2016.

\end{thebibliography}

	\begin{appendix}
		\phantomsection
		\setcounter{dummythm}{0}%
		\setcounter{equation}{0}%
		\renewcommand{\theHdummythm}{A.\arabic{dummythm}}%
		\renewcommand{\theHequation}{A.\arabic{equation}}%
		\renewcommand{\thedummythm}{A.\arabic{dummythm}}%
		\renewcommand{\theequation}{A.\arabic{equation}}%
\begin{lem}\label{thm:strictdiff}%
	Suppose that \(\func{G}{\R^n}{\R^p}\) is locally Lipschitz around a point \(\bar z\) at which \(G(\bar z)=0\).
	Then, \(\phi(z)\coloneqq\tfrac12\|G(z)\|^2\) is strictly differentiable at \(\bar z\) (in the sense of \cite[Def. 9.17]{rockafellar2009variational}) with null gradient.
	\begin{proof}
		Let \(L\) be a Lipschitz constant for \(G\) in a neighborhood \(\mathcal U\) of \(\bar z\).
		Then, for \(z,z'\in\mathcal U\) we have
		\begin{align*}
		&
			\frac{\bigl|
				\phi(z)-\phi(z')-\innprod{0}{z-z'}
			\bigr|}{\|z-z'\|}
		{}={}
			\frac{\bigl|
				\|G(z)\|^2-\|G(z')\|^2
			\bigr|}{2\|z-z'\|}
		\\
		{}={} &
			\frac{\bigl|
				\|G(z)-G(z')\|^2+2\innprod{G(z')}{G(z)-G(z')}
			\bigr|}{2\|z-z'\|}
		\\
		{}\leq{} &
			\tfrac L2\|G(z)-G(z')\|
			{}+{}
			L\|G(z')\|
		\end{align*}
		which vanishes as \(z,z'\to\bar z\), hence the claim.
	\end{proof}
\end{lem}

\begin{lem}\label{lem:quadratic}
	Let $\mathcal H_i\in\R^{r_i\times r_i}$ be symmetric positive definite, $\mathcal V_i\in\R^{r_i\times p}$, and $\nu_i\in\R^{r_i}$, $i\in[N]$.
	If $\mathcal U\coloneqq\sum_{i=1}^N\mathcal V_i^{\top}\mathcal H_i^{-1}\mathcal V_i$ is symmetric positive definite, then
	\[\textstyle
		\sum_{i=1}^N\|\mathcal V_ix-\nu_i\|_{\mathcal H_i^{-1}}^2
	{}={}
		\|\mathcal Ux-d\|_{\mathcal U^{-1}}^2
		{}-{}
		\|d\|_{\mathcal U^{-1}}^2
		{}+{}
		\sum_{i=1}^N\|\nu_i\|_{\mathcal H_i^{-1}}^2,
	\]
	where $d\coloneqq\sum_{i=1}^N\mathcal V_i^{\top}\mathcal H_i^{-1}\nu_i$.
	\begin{proof}
		Let $q(x)= \sum_{i=1}^N\|\mathcal V_ix-\nu_i\|_{\mathcal H_i^{-1}}^2$.
		That
		\(
			q(x)
		{}={}
			\|x\|_{\mathcal U}^2 + \sum_{i=1}^N\|\nu_i\|_{\mathcal H_i^{-1}}^2 - 2\innprod xd
		\)
		is of immediate verification.
		Since $q$ is quadratic, the Taylor expansion around its minimizer $x^\star = {\mathcal U}^{-1}d$ is given by $q(x) = q(x^\star) + \|x - x^\star\|_{\mathcal U}^2$.
		Substituting $x^\star$ results in the claimed form.
	\end{proof}
\end{lem}

\begin{appendixproof}{thm:ALM}
	Owing to the update at \cref{state:ALM:y},
\[
	\LL(\bm z,\bm\lambda^{k+1})
{}={}
	\LL_{\beta_k}(\bm z,\bm\lambda^k)
	{}+{}
	\tfrac{\beta_k}2\|F(\bm z^{k+1})\|^2
	{}-{}
	\tfrac{\beta_k}2\|F(\bm z)-F(\bm z^{k+1})\|^2.
\]
The last term on the right-hand side is continuously differentiable (with null gradient) at \(\bm z^{k+1}\), owing to \cref{thm:strictdiff}.
It then follows from \cite[Ex. 8.8(c)]{rockafellar2009variational} that
\(
	\partial_{\bm z}\LL(\bm z^{k+1},\bm\lambda^{k+1})
{}={}
	\partial_{\bm z}\LL_{\beta_k}(\bm z^{k+1},\bm\lambda^k)
\),
hence that the pair \((\bm z^k,\bm\lambda^k)\) satisfies condition \eqref{eq:KKTopt} for every \(k\geq1\), by virtue of \cref{prop:ALM:grad} in \cref{state:ALM:x}.
It remains to show that \eqref{eq:KKTfeas} too is eventually satisfied.
Notice that, by definition of \(\hat z^k\) at \cref{state:ALM:barx}, \(\LL(\bm z^{k+1},\bm\lambda^k)\leq f(\bm z^0)\) holds for every \(k\), which combined with \eqref{eq:L'} yields
{\mathtight[0.85]\begin{align*}
	\tfrac{1}{2\beta_k}
	\|\bm\lambda^{k+1}\|^2
{}={} &
	\tfrac{\beta_k}{2}
	\|F(\bm z^{k+1})+\bm\lambda^k\nicefrac{}{\beta_k}\|^2
{}\leq{}
	f(\bm z^0)
	{}-{}
	f(\bm z^{k+1})
	{}+{}
	\tfrac{1}{2\beta_k}\|\bm\lambda^k\|^2
\\
{}\leq{} &
	c
	{}+{}
	\tfrac{1}{2\beta_k}\|\bm\lambda^k\|^2,
\end{align*}}%
where \(c\coloneqq f(\bm z^0)-\inf f\) is a constant.
\begin{subequations}\label{subeq:dual_bound}%
	Since \(\beta_{k+1}\geq\beta_k\), it holds that
	\(
	\tfrac{1}{2\beta_{k+1}}
		\|\bm\lambda^{k+1}\|^2
		{} \leq {}
		c +
		\tfrac{1}{2\beta_k}
		\|\bm\lambda^k\|^2
	\),
	which leads to
	\begin{equation}
		\tfrac{1}{\beta_k}\|\bm\lambda^k\|^2
	{}\leq{}
		\tfrac{1}{\beta_0}\|\bm\lambda^0\|^2
		{}+{}
		2kc
	\end{equation}
	for every \(k\in\N\).
	Moreover, since
	\begin{align}
	\nonumber
		\tfrac12\|F(\bm z^{k+1})\|^2
	{}\leq{} &
		\|F(\bm z^{k+1})+\bm\lambda^k\nicefrac{}{\beta_k}\|^2
	{}+{}
		\|\bm\lambda^k\nicefrac{}{\beta_k}\|^2
	\\
	{}={} &
		\tfrac{1}{\beta_k^2}\bigl[
			\|\bm\lambda^{k+1}\|^2
	{}+{}
			\|\bm\lambda^k\|^2
		\bigr],
	\end{align}
	the \(\beta\)-update at \cref{state:ALM:beta} implies that \(\|F(\bm z^{k+1})\|_\infty\to0\) (\(Q\)-linearly) if \(\beta_k\) is asymptotically constant, hence the claim.
\end{subequations}
Otherwise, the set \(\mathcal K\coloneqq\set{k\in\N}[\beta_k=\max\set{\xi\beta_{k-1},\beta_0k^\alpha}]\) is infinite.
Then, for \(k\in\mathcal K\), combining \eqref{subeq:dual_bound} yields
\begin{align*}
	\tfrac12\|F(\bm z^{k+1})\|^2
{}\leq{} &
	\textstyle
	\tfrac{\beta_{k+1}}{\beta_k^2}\bigl(\tfrac{1}{\beta_0}\|\bm\lambda^0\|^2+2(k+1)c\bigr)+\tfrac{1}{\beta_k}\bigl(\tfrac{1}{\beta_0}\|\bm\lambda^0\|^2+2kc\bigr)
\\
{}\leq{} &
	\max\set{\tfrac{\xi}{\beta_k},\tfrac{\beta_0(k+1)^\alpha}{\beta_k^2}}\bigl(\tfrac{1}{\beta_0}\|\bm\lambda^0\|^2 + 2(k+1)c\bigr)
\\
&
	{}+{}
	\tfrac{1}{\beta_k}\bigl(\tfrac{1}{\beta_0}\|\bm\lambda^0\|^2+2kc\bigr)
\\
	\dueto{(\(\beta_k\geq k^\alpha\))}
	~
{}\leq{} &
	\max\set{\tfrac{\xi}{k^\alpha},\tfrac{\beta_0(k+1)^\alpha}{k^{2\alpha}}}\bigl(\tfrac{1}{\beta_0}\|\bm\lambda^0\|^2 + 2(k+1)c\bigr)
\\
&
	{}+{}
	\tfrac{1}{k^\alpha}\bigl(\tfrac{1}{\beta_0}\|\bm\lambda^0\|^2+2kc\bigr)
{}\to{}
	0
\end{align*}
as \(\mathcal K\ni k\to\infty\), owing to the fact that \(\alpha>1\).
The second inequality uses the fact that, regardless of whether $k+1\in\mathcal K$ or not, $\beta_{k+1}\leq \max\set{\xi\beta_k,\beta_0(k+1)^\alpha}$ holds (since $\xi>1$).

\end{appendixproof}

\begin{appendixproof}{thm:GN:ST}
Note that matrices $A_{\Ell}$, $B_{\Ell}$ and vectors $\bs_j$ in \cref{alg:GN} depend on the current iterate $\bm z^l$.
Here, we use superscript $l$ to emphasize this dependance.
The linear least squares \cref{prob:GN} solved at \cref{state:GN:min} may equivalently be written as
\begin{equation}\label{eq:LSform}
	\minimize_{\bm z = (\bm w, \bm x)}{
		\tfrac12\|J(\bm z^l)\bm z - b(\bm z^l)\|^2,
	}
\end{equation}
where
\begin{align*}
	J(\bm z^l)
{}={} &
	\begin{pmatrix}
		-\bar B^l & \bar A^l \\
		\sqrt{\mu_w}\I &
	\end{pmatrix}
\shortintertext{%
	with
}
	\bar A
{}={} &
	\renewcommand{\arraystretch}{1.2}
	\begin{pmatrix}
		\sqrt{\rho_1}\I																		\\
		-\sqrt{\rho_2}A_2^l	& \sqrt{\rho_2}\I 												\\
							& \ddots			& \ddots									\\
							&					& \mathllap{-{}}\sqrt{\rho_{N}}A_N^l	& \hspace{-25pt}\sqrt{\rho_{N}}\I	\\
							&					&						& \hspace{-15pt}-\sqrt{\rho_{N+1}}A_{N+1}^l
	\end{pmatrix},
\\
	\bar B^l
{}={} &
	\blkdiag\bigl(\sqrt{\rho_1}B_1^l,\dots, \sqrt{\rho_{N+1}}B_{N+1}^l\bigr),
\end{align*}
and $b(\bm z^l) = (\sqrt{\rho_1}c_1,\ldots, \sqrt{\rho_{N+1}}c_{N+1}, 0_s \big)$.

In what follows we show that the eigenvalues of $J(\bm z^l)^\top J(\bm z^l)$ along a converging subsequence are bounded above and away from zero, and that \cref{state:GN:LS} is a restatement of the standard Armijo linesearch, at which point the claim follows from standard results for gradient methods \cite{bertsekas2016nonlinear,facchinei2003finite}.
For the latter, note that by the optimality conditions for \eqref{eq:LSform}, the solution $\bar{\bm z}^l$ satisfies
\(
	J(\bm z^l)^\top J(\bm z^l) \bm{\bar{z}}^l = J(\bm z^l)^\top b(\bm z^l).
\)
Moreover, since $\nabla_z\LL_\beta(\bm z^l,\bl)=J(\bm z^l)^{\top}J(\bm z^l) \bm z^l - J(\bm z^l)^{\top} b(\bm z^l)$, combining the two equalities yields
\(
	\innprod{\nabla_z\LL_\beta(\bm z^l,\bl)}{\bar{\bm z}^l - \bm z^l}
{}={}
	-\|J(\bm z^l)(\bar{\bm z}^l - \bm z^l)\|^2,
\)
establishing the claimed equivalence.

Let $\seq{\bm z^l}[l\in K]$ be a subsequence converging to a limit point $\bm z^\star$.
Note that $J(\bm z)$ has full column rank for any $\bm z$, and $J(\bm z)^{\top} J(\bm z)$ is thus nonsingular.
Therefore, by continuity of $J(\cdot)$ and \cite[Lem. 7.5.2]{facchinei2003finite} we have that $c_1\I \preceq J(\bm z^l)^{\top} J(\bm z^l) \preceq c_2 \I$ for some $c_1, c_2>0$.
The claim then follows from \cite[Prop. 8.3.7]{facchinei2003finite}.

\end{appendixproof}

\begin{appendixproof}{thm:FDP}
First, note that by \eqref{eq:E_G}
\begin{align*}
	G_{\Ell}^{\top}G_{\Ell}+\tfrac{\mu_w}{\rho_{\Ell}} E_{\Ell}^{\top}E_{\Ell}
{}={} &
	\I-B_{\Ell}E_{\Ell}-E_{\Ell}^{\top}B_{\Ell}^{\top}+E_{\Ell}^{\top}\left(B_{\Ell}^{\top}B_{\Ell}+\tfrac{\mu_w}{\rho_{\Ell}}\I\right)E_{\Ell}
\\
\numberthis\label{eq:M_D}
{}={} &
	\I-B_{\Ell}E_{\Ell} {}={} G_{\Ell}.
\end{align*}
We proceed by induction to show that, for $\Ell\in[N]$,
\begin{equation}\label{eq:Vk:induction}
	V_\Ell^{\star}(\bx_{\Ell}) = \tfrac12\|M_{\Ell}^{-1}\bx_{\Ell}-\bq_{\Ell}\|_{M_{\Ell}}^2 + C_{\Ell},
\end{equation}
where the term $C_{\Ell}$ does not depend on $\bx_{\Ell}, \bx_{\Ell+1},\ldots,\bx_N$.
Here, we avoid deriving a recursion for $C_{\Ell}$ since it does not affect the computation of $\bx_{\Ell}$ and $\bw_{\Ell+1}$ in the next stages.

For the base case $\Ell=1$, by the first order optimality condition for the minimization \eqref{eq:FDR:V1} the unique minimizer is computed as
\[
	\bw_1^\star(x_1) = E_1 \big(\bx_1-A_1\bx_{0}-\bs_1\big).
\]
After substitution, using \eqref{eq:M_D}, and simple algebra we obtain
\[
	V_1^{\star}(\bx_1) {}={} \tfrac12\|\bx_1-A_1\bx_{0}-\bs_1\|_{\rho_1G_1}^2 {}={}\tfrac12\|M_1^{-1}\bx_1-\bq_1\|_{M_1}^2
\]
where $M_1=\rho_1^{-1}G_1^{-1}$, and $q_1 = M_1^{-1}(A_1\bx_{0}+\bs_1)$.

Arguing by induction, suppose that \eqref{eq:Vk:induction} holds for some $\Ell$ such that $1\leq \Ell\leq N-1$.
Let $\varphi(\bx_{\Ell},\bw_{\Ell+1})$ denote the argument being minimized in \eqref{eq:FDR:Vk}.
From direct computation
\[
	\nabla^2 \varphi(\bx_{\Ell},\bw_{\Ell+1})
{}={}
	\begin{pmatrix}
		M_{\Ell}^{-1}+\rho_{\Ell+1}A_{\Ell+1}^{\top}A_{\Ell+1} & \rho_{\Ell+1}A_{\Ell+1}^{\top}B_{\Ell+1}
	\\
		\rho_{\Ell+1}B_{\Ell+1}^{\top}A_{\Ell+1} & \mu_w\I+\rho_{\Ell+1}B_{\Ell+1}^{\top}B_{\Ell+1}
	\end{pmatrix}.
\]
Since $M_\Ell\in\sym{r_\Ell}_{++}$, by forming its Schur complement and using \eqref{eq:E_G} it follows that the Hessian is symmetric positive definite if and only if so is $M_{\Ell}^{-1}+\rho_{\Ell+1}A_{\Ell+1}^{\top}G_{\Ell+1}A_{\Ell+1}$, which holds true.
Hence, the subproblems have unique solutions.
By the first order optimality condition for \eqref{eq:FDR:Vk}, the solution pair $(\bx_{\Ell}^\star, \bw^\star_{\Ell+1})$ satisfies
{\mathtight[0.9]\begin{align*}
	0
{}={} &
	M_{\Ell}^{-1}\bx_{\Ell}^\star-\bq_{\Ell}
	{}-{}
	\rho_{\Ell+1}A_{\Ell+1}^{\top}
	\left(
		\bx_{\Ell+1}-A_{\Ell+1}\bx_{\Ell}^\star-B_{\Ell+1}\bw_{\Ell+1}^\star-\bs_{\Ell+1}
	\right)
\shortintertext{%
	and
}
	0
{}={} &
	\mu_w \bw_{\Ell+1}^\star-\rho_{\Ell+1} B_{\Ell+1}^{\top}\left(\bx_{\Ell+1}-A_{\Ell+1}\bx_{\Ell}^\star-B_{\Ell+1}\bw_{\Ell+1}^\star-\bs_{\Ell+1}\right).
\end{align*}}%
The latter reads $\bw_{\Ell+1}^\star = E_{\Ell+1}\big(\bx_{\Ell+1}-A_{\Ell+1}\bx_{\Ell}^\star-\bs_{\Ell+1}\big)$.
After substituting $\bw_{\Ell+1}^\star$ into the former, using \eqref{eq:E_G} and \eqref{eq:M_D} we obtain
\[
	\bx_{\Ell}^{\star} = S_{\Ell+1}\bq_{\Ell}+P_{\Ell}(\bx_{\Ell+1}-\bs_{\Ell+1}),
\]
where $P_k=\rho_{\Ell+1}S_{\Ell+1}A_{\Ell+1}^{\top}G_{\Ell+1}$ and $S_{\Ell+1}$ is as in \eqref{eq:MatrixUpdates_S}.
Substituting the minimizer pair $(\bx_{\Ell}^\star, \bw_{\Ell+1}^\star)$ back in \eqref{eq:FDR:Vk} and using \eqref{eq:M_D} yields
\begin{align*}
	V_{\Ell+1}^\star(\bx_{\Ell+1})
{}={} &
	C_{\Ell} + \tfrac12\|\M_{\Ell}^{-1}\bx_{\Ell}^{\star}-\bq_{\Ell}\|_{\M_{\Ell}}^2 \\
&
	+
	\tfrac12\|\bx_{\Ell+1}-A_{\Ell+1}\bx_{\Ell}^{\star}-\bs_{\Ell+1}\|_{\rho_{\Ell+1}G_{\Ell+1}}^2 \\
	\dueto{(subs. $\bx_{\Ell}^\star$)}
		{}={} &
	C_{\Ell} + \tfrac12\|\mathcal{V}_1\bx_{\Ell+1}-\nu_1\|_{\mathcal{H}_1^{-1}}^2 \\
		&{}+{}
	\tfrac12\|\mathcal{V}_{2}\bx_{\Ell+1}-\nu_{2}\|_{\mathcal{H}_{2}^{-1}}^2
\end{align*}
where $\mathcal{V}_1=M_{\Ell}^{-1}P_{\Ell}$, $\mathcal{H}_1=M_{\Ell}^{-1}$, $\mathcal{V}_{2}=\left(\I-A_{\Ell+1}P_{\Ell}\right)$, $\mathcal H_{2}=\rho_{\Ell+1}^{-1}G_{\Ell+1}^{-1}$,
\begin{align*}
	\nu_1
{}={} &
	\left(\I-M_{\Ell}^{-1}S_{\Ell+1}\right)\bq_{\Ell}+M_{\Ell}^{-1}P_{\Ell}\bs_{\Ell+1},
	~\text{and}
\\
	\nu_{2}
		{}={} &
	A_{\Ell+1}S_{\Ell+1}\bq_{\Ell}+(\I-A_{\Ell+1}P_{\Ell})\bs_{\Ell+1}.
\end{align*}
On the other hand, we have that
\begin{align*}
	\mathcal{U}
		{}={} &
		\sum_{i=1}^2\mathcal{V}_i^{\top}\mathcal{H}_i^{-1}\mathcal{V}_i \\
		{}={} &
		P_{\Ell}^{\top}\M_{\Ell}^{-1}P_{\Ell}+\rho_{\Ell+1}\left(\I-A_{\Ell+1}P_{\Ell}\right)^{\top}G_{\Ell+1}\left(\I-A_{\Ell+1}P_{\Ell}\right) \numberthis \label{eq:Malt} \\
		{} = {} &
		P_{\Ell}^{\top}S_{\Ell+1}^{-1}P_{\Ell}+\rho_{\Ell+1}\big(\I-P_{\Ell}^{\top}A_{\Ell+1}^{\top}\big)G_{\Ell+1}-\rho_{\Ell+1}G_{\Ell+1}A_{\Ell+1}P_{\Ell} \\
		{}={} &
		\rho_{\Ell+1}G_{\Ell+1}-\rho_{\Ell+1}^2G_{\Ell+1}A_{\Ell+1}S_{\Ell+1}A_{\Ell+1}^{\top}G_{\Ell+1}
		\\
		{} = {} &
		M_{\Ell+1}^{-1},
\end{align*}
where \eqref{eq:SnM} was used in the second equality, and the Woodbury matrix identity was used in the last equality.
Therefore, we may apply \cref{lem:quadratic} to obtain
\begin{align}
\nonumber
	V_{\Ell+1}^\star(\bx_{\Ell+1})
{}={} &
	\tfrac12\|M_{\Ell+1}^{-1}\bx_{\Ell+1}-\bq_{\Ell+1}\|_{M_{\Ell+1}}^2
\\
\label{eq:Vk+}
&
	\textstyle
	{}-{} \tfrac12\|\bq_{\Ell+1}\|_{M_{\Ell+1}}^2+\tfrac12\sum_{i=1}^2\|\nu_i\|_{\mathcal{H}_i^{-1}}^2 + C_{\Ell},
\end{align}
with
\begin{align}
	\bq_{\Ell+1}
	{}={} &
	P_{\Ell}^{\top}\nu_1+\rho_{\Ell+1}\left(\I-A_{\Ell+1}P_{\Ell}\right)^{\top}G_{\Ell+1}\nu_{2}
	\\
	{}={} &
	M_{\Ell+1}^{-1}\bs_{\Ell+1}+\rho_{\Ell+1}G_{\Ell+1}A_{\Ell+1}\bq_{\Ell},
\end{align}
where we used \eqref{eq:SnM} and the alternative expression for $M_{\Ell+1}$ in \eqref{eq:Malt}.
The last three terms in \eqref{eq:Vk+} are absorbed into $C_{\Ell+1}$ completing the induction argument.
A recursive formula for $C_{\Ell}$ is not provided since it does not depend on future states and as such would not effect the solution to the minimization of the next stages.

It remains to solve \eqref{eq:FDR:Vfinal}.
Arguing as before, from the first order optimality condition the solution pair $(\bx_N^\star, \bw^\star_{N+1})$ must satisfy
\begin{align*}
	0
		{}={} &
	M_N^{-1}\bx_N^\star-\bq_N \\
	& + \rho_{N+1}A_{N+1}^{\top}\left(A_{N+1}\bx_N^{\star}+B_{N+1}\bw_{N+1}^{\star}+\bs_{N+1}\right),
\end{align*}
and
\[
	0 = \mu_w \bw_{N+1}^\star + \rho_{N+1} B_{N+1}^{\top}\left(A_{N+1}\bx_N^{\star}+B_{N+1}\bw_{N+1}^{\star}+\bs_{N+1}\right),
\]
The former equality is equivalent to the one given in \cref{state:Backward}.
After substituting $\bw_{N+1}^\star$ back into the latter and using \eqref{eq:M_D}, the update for $\bx_N^\star$ is obtained.

\end{appendixproof}

	\end{appendix}

\end{document}